%% file: hyprsfc.tex
\documentclass{article}

\input{header}

\newcommand{\A}{{\mathcal A}}
\newcommand{\proj}[3]{\pi_{#3}(#2,#1)}

\title{Sharp $L^p$-$L^q$ estimates for generalized $k$-plane transforms}

\author{Philip T. Gressman}
\begin{document}
\maketitle
\begin{abstract}
In this paper, optimal $L^p-L^q$ estimates are obtained for operators which average functions over polynomial submanifolds, generalizing the $k$-plane transform.  An important advance over previous work (e.g., \cite{gressman2006}) is that full $L^p-L^q$ estimates are obtained by methods which have traditionally yielded only restricted weak-type estimates.  In the process, one is lead to make coercivity estimates for certain functionals on $L^p$ for $p < 1$.
\end{abstract}
\section{Introduction}
\subsection{Background}

The object of study in this paper is the family of operators which integrate a function $f$ over all submanifolds given by polynomials in some appropriate coordinate system.  On $\R^2$, for example, such an operator would map a polynomial $p$ in a single variable to the integral of a function $f(x,y)$ over the graph $y = p(x)$.  To be more precise, 
fix any positive integers $n$, $n'$, and $d$. Let $M_{n,d}$ be the set of all multiindices of length $n$ and degree at most $d$ (recall that a multiindex is simply an $n$-tuple of nonnegative integers, the degree of a multiindex $\alpha = (\alpha_1,\ldots,\alpha_n)$ is denoted by $|\alpha|$ and equals $\sum_{j=1}^n \alpha_j$, and $t^\alpha := \prod_{j=1}^n t_j^{\alpha_j}$ as well as $\alpha! := \prod_{j=1}^n \alpha_j!$).  Let $T_{n,n',d}$ be the operator mapping functions on $\R^{n} \times \R^{n'}$ (written $f(x,y)$ where $x \in \R^n$ and $y \in \R^{n'}$) to functions on $(\R^{n'})^{M_{n,d}}$ (thought of as the space of coeffiecients of an $n'$-tuple of polynomials in $n$ variables of degree at most $d$) given by\begin{equation}
T_{n,n',d} f(u) := \int_{\R^n} f \left( t, \sum_{\alpha \in M_{n,d}} u_\alpha t^\alpha \right) dt  \label{theop}
\end{equation}
(i.e., $u_\alpha \in \R^{n'}$ for all $\alpha \in M_{n,d}$ and $u=(u_\alpha)_{\alpha \in M_{n,d}}$).  The purpose of this paper is to establish the $L^p$ mapping properties of the family \eqref{theop}.

This family \eqref{theop}, when $d=1$, generalizes the classical $k$-plane transform.  %The operator $T_{1,n-1,1}$ corresponds to the classical X-ray transform on $\R^n$.  Similarly, $T_{k,n-k,1}$ is essentially the $k$-plane transform on $\R^n$.  
To see this, let $u_0, \ldots, u_k \in \R^{n-k}$ be vectors, and consider the following mapping into ${\cal M}_{k,n}$, the space of all affine $k$-planes in $\R^n$:
\[\sigma(u_{0}, u_{1}, \ldots, u_{k}) := \set{(t,u_0 + u_1 t_1 + \cdots + u_k t_k) \in \R^n}{t = (t_1,\ldots,t_k) \in \R^k}. \]
Provided that $||u_1|| + \cdots + ||u_k||$ is small, the pull-back of the natural measure on ${\cal M}_{k,n}$ is comparable to the Lebesgue measure on $\R^{(n-k)(k+1)}$; furthermore, the the pull-back of the Lebesgue measure on the $k$-plane $\sigma(u)$ is comparable to $dt$, that is, 
\[ C^{-1} T_{k,n-k,1} f(u) \leq \int_{\sigma(u)} f \leq  C T_{k,n-k,1}f(u)\]
(where $f$ is a nonnegative function).  %As long as $||u_1|| + \cdots + ||u_k||$ is small, this inequality may be integrated as well.  
To obtain global inequalities, one simply averages over all rotations $f^\theta(x) := f(\theta \cdot x)$ for $\theta \in SO(n)$.

The $L^p$-boundedness (in both standard and mixed-norm spaces) of the classical $k$-plane transform (including the Radon transform as a special case) was established in the 1980s in various papers including, for example, the works of Christ \cite{christ1984}, Drury \cite{drury1983}, \cite{drury1984}, \cite{drury1986}, and Oberlin and Stein \cite{os1982}.  The classical estimate to be generalized by this paper is that the $k$-plane transform maps $L^p(\R^n)$ to $L^q({\cal M}_{k,n})$ when $p= \frac{n+1}{k+1}$ and $q = n+1$ (the restricted weak-type version was established by Drury \cite{drury1983}, and the full estimate by Christ \cite{christ1984}).  By the remarks of the preceding paragraph, this estimate is a special case of theorem \ref{maintheorem} after performing the prescribed average over rotations. % are a consequence the main theorem of this paper after performing the prescribed average over rotations.  %, and $T_{k,n-k,1}$ is simply a weighted version of the $k$-plane transform on $\R^{n}$.

When $d > 1$, the corresponding operators are largely new.  The family $T_{1,n-1,d}$ arose in earlier work of the author \cite{gressman2006} as examples of overdetermined, one-dimensional averaging operators.  The significance of the family $T_{1,n-1,d}$ in that paper is that such operators are, in some sense, less degenerate than the classical X-ray transform; for example, the operators $T_{1,n-1,d}$ map $L^p_{comp}$ to $L^q_{loc}$ for a larger set of indices $(p,q)$ than does the classical X-ray transform.  The main result of that paper \cite{gressman2006} was a family of restricted weak-type estimates; in this paper, the corresponding strong estimates follow from theorem \ref{maintheorem} as well. % those estimates are replaced by full $L^p-L^q$ estimates.
%Modulo certain restricted weak-type inequalities, the $L^p$ theory of these operators was settled in that earlier work.

The Fourier integral operator realization of \eqref{theop} has nondegenerate canonical relation, so earlier theorems concerning overdetermined averaging operators, including recent work of Brandolini, Greenleaf, and Travaglini \cite{bgt2005} and Ricci and Travaglini \cite{rt2001}, can be applied.  The proofs of these theorems are heavily concerned with the behavior of the operators \eqref{theop} near $L^2$ (and rely on oscillatory integral estimates in one form or another).  Such theorems give suboptimal results in general, meaning that they are restricted to the study of $L^p \rightarrow L^{p'}$ estimates for conjugate exponents $p,p'$.  Unlike these earlier results, this paper approaches the question from the standpoint of geometric combinatorics (pioneered by Christ \cite{christ1998}), and is able to establish complete results. % concerning the $L^p$ regularity of the operators \eqref{theop}.

It is useful to note that the operators \eqref{theop} possess a variety of symmetries.  First and foremost is an $(n + n')$-dimensional family of dilation symmetries:  taking $f_{\delta,\delta'}(x,y) := f(\delta_1 x_1, \ldots, \delta_n x_n, \delta'_1 y_1,\ldots, \delta'_{n'} y_{n'})$ for arbitrary positive numbers $\delta_1,\ldots,\delta_n,\delta'_1,\ldots,\delta'_{n'}$ induces a scaling $u^{\delta,\delta'}$ on the space $(\R^{n'})^{M_{n,d}}$ by requiring $T_{n,n',d} f_{\delta,\delta'} (u) = T_{n,n',d} f (u^{\delta,\delta'})$.  This family will appear explicitly in section \ref{necessitysec}.  Likewise, the translations $f_{h,h'}(x,y) := f(x+h,y+h')$ induce a family of translation operators $\tau^{h,h'}$ on $(\R^{n'})^{M_{n,d}}$ by (again) requiring that $T_{n,n',d} f_{h,h'} (u) = T_{n,n',d} f(\tau^{h,h'}(u))$.  Although this family $\tau^{h,h'}$ is not the usual family of translations on Euclidean space, it does possess many of the same properties, including that $\tau^{h,h'}$ is measure-preserving with respect to the Lebesgue measure $d u$.  Finally, functions in the range of $T_{n,n',d}$ satisfy a family of PDEs (a fact first observed by F. John \cite{john1938}).  Let $\partial_{\alpha,j}$ be differentiation with respect to the $j$-th component of $u_{\alpha}$ for $\alpha \in M_{n,d}$.  For any $j$ and $k$ between $1$ and $n'$ (inclusive) and any $\alpha, \beta, \tilde \alpha, \tilde \beta \in M_{n,d}$ satisfying $\alpha + \beta = \tilde \alpha + \tilde \beta$,
\[\left[\partial_{\alpha,j} \partial_{\beta,k} - \partial_{\tilde{\alpha},j} \partial_{\tilde{\beta},k} \right] T_{n,n',d} f (u) = 0 \]
(in the sense of distributions) for any $f$.
\subsection{Main theorems}
\label{mainthms}
The main theorems of this paper establish sharp $L^p-L^q$ boundedness of \eqref{theop} and related generalizations. As already mentioned, an important feature (not found previously) of these theorems is that full endpoint estimates are obtained, not simply restricted weak-type estimates.
The first theorem deals with the global and local $L^p-L^q$ mapping properties of the family \eqref{theop}: 

\begin{theorem}
The operator $T_{n,n',d}$ maps $L^p(\R^{n} \times \R^{n'})$ to $L^q ((\R^{n'})^{M_{n,d}})$ if and only if $p = 1 + \frac{n' d}{n+1}$ and $q = |M_{n,d}| p$, where $|M_{n,d}| = \binom{n+d}{d} := \frac{(n+d)!}{n! d!}$ is the number of multiindices of length $n$ and degree at most $d$.  Furthermore, $T_{n,n',d}$ maps $L^p_{comp} \rightarrow L^q_{loc}$ if and only if $(p^{-1},q^{-1})$ is in the closed convex hull of the points $(0,0)$, $(1,1)$, $(0,1)$ and 
\[ \left( \frac{n+1}{n + n' j + 1}, { \binom{n+j}{j} }^{-1} \frac{(n+1)}{ (n + n'j + 1)} \right) \] \label{maintheorem}
for $j=1,\ldots,d$.
\end{theorem}

\begin{figure}
\centering
\psset{unit=6cm}
\begin{pspicture}(-0.1,-0.1)(1.1,1.1)
%\put(3.1,1.4){$\scriptstyle (p_1^{-1},1-{q'_1}^{-1})$}
%\put(1.3,0.1){ $ (p_j^{-1},1-{q'_j}^{-1})$}
%\put(1.6,0.3){$\scriptstyle \cdots$}
%\put(0,0){\line (1,1){200}}
\pspolygon*[linecolor=lightgray](0,0)(0,1)(1,1)(0.75,0.25)(0.6,0.1)(0.5,0.05)(0.42857,0.02857)(0.375,0.01785)(0,0)
\pspolygon[linewidth=0.5pt](0,0)(0,1)(1,1)(0.75,0.25)(0.6,0.1)(0.5,0.05)(0.42857,0.02857)(0.375,0.01785)(0,0)
\psline[linewidth=0.5pt](0,0)(1,0)(1,1)(0,1)(0,0)
%\psline[linewidth=0.5pt](0.5,0)(1,0.5)(0.5,1)(0,0.5)
\put(-0.08,0.48){$\displaystyle \frac{1}{q}$}
\put(0.45,-0.075){$\displaystyle 1/p$}
%\put(-0.07,-0.06){$\displaystyle (0,0)$}
\put(-0.005,-0.06){$\displaystyle 0$}
\put(0.98,-0.06){$\displaystyle 1$}
\put(-0.04,0.96){$\displaystyle 1$}
\put(-0.04,0.01){$\displaystyle 0$}
\pscircle*(0.75,0.25){0.005}
\pscircle*(0.6,0.1){0.005}
\pscircle*(0.5,0.05){0.005}
\pscircle*(0.42857,0.02857){0.005}
\pscircle*(0.375,0.01785){0.005}
%\pscircle*(1,5.833){0.05}
%\pscircle*(0.857,5.878){0.05}
\end{pspicture}
\caption{The operator $T_{2,1,5}$ maps $L^p_{comp}$ to $L^q_{loc}$ precisely when $(p^{-1},q^{-1})$ is in the shaded polygon shown above. The nontrivial vertices are marked by dots.}
\label{thefig}
\end{figure}
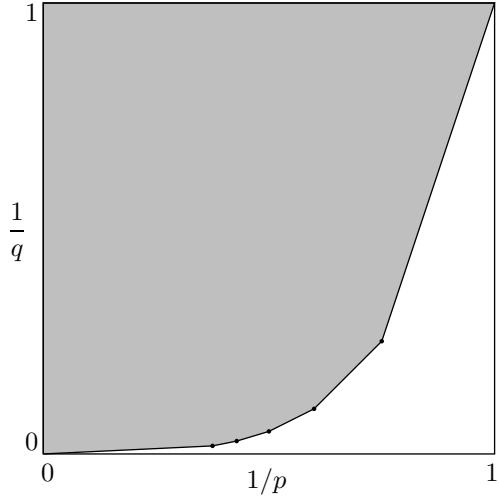

As indicated by the theorem (see figure \ref{thefig} for an illustration of the Riesz diagram of a typical operator), the local mapping properties of $T_{n,n',d}$ are far more complex than a non-overdetermined averaging operator with nonvanishing rotational curvature (for example).  If one moves to the scale of mixed-norm spaces, the boundedness properties of $T_{n,n',d}$ become even more complex because of natural ``factorizations'' which occur.  For example,
to prove the $L^p_{comp} \rightarrow L^q_{loc}$ estimates, it suffices to prove {\it only} global estimates.  This is because an $L^p \rightarrow L^q$ estimate for $T_{n,n',j}$ implies mixed-norm boundedness of the form $L^p \rightarrow L^\infty(L^q)$ for $T_{n,n',d}$ when $d > j$, where the $L^q$-norm is taken over the variables $u_\alpha$ for $|\alpha| \leq j$.  To see this, let $\pi$ be the natural projection from $(\R^{n'})^{M_{n,d}}$ to $(\R^{n'})^{M_{n,j}}$ and let $\widetilde{\pi}$ be the corresponding projection onto the orthogonal complement of $(\R^{n'})^{M_{n,j}}$.  It follows that
\begin{equation}
 T_{n,n',d} f(u) = T_{n,n',j} f^{\widetilde{\pi}(u)} (\pi(u))\label{split}
\end{equation}
where $f^{\widetilde{\pi}(u)}(t,x) := f(t,x + \sum_{|\alpha| > j} u_\alpha t^\alpha)$.
Thus
\[ \int_{\widetilde{\pi}(u) = v} |T_{n,n',d} f(u)|^q \leq C ||f^v||^q_p = C ||f||_p^q \]
uniformly in $v$ when $T_{n,n',j}$ maps $L^p \rightarrow L^q$.

To identify many of the ``trivial'' $L^\infty(L^q)$ estimates satisfied by the operators \eqref{theop}, it is useful to restrict $T_{n,n',d}$ to somewhat general hyperplanes.  In this paper, attention is fixed on coordinate hyperplanes.  Such hyperplanes will be identified by the coordinate axes they contain; the axes themselves will be identified with elements in $M_{n,d} \times \{1,\ldots,n'\}$, so coordinate hyperplanes are identified with subsets $\A \subset M_{n,d} \times \{1,\ldots,n'\}$.
For each $j=1,\ldots,n'$, let $\A_j$ be the collection of multiindices $\alpha$ for which $(\alpha,j) \in \A$.  Then the restriction of $T_{n,n',d}$ to the coordinate hyperplane given by $\A$ will be denoted $T_\A$; the explicit formula is simply %there is an operator $T_\A$ associated to $\A$ (which is a restriction of $T_{n,n',d}$ to a coordinate hyperplane), given by
\begin{equation}
T_\A f(u) := \int_{\R^n} f \left( t, \sum_{\alpha \in \A_1} u_{(\alpha,1)} t^\alpha, \ldots, \sum_{\alpha \in \A_{n'}} u_{(\alpha,n')} t^\alpha \right) dt; \label{theop2}
\end{equation}
for convenience, the following shorthand will be used in the future:
\[
\proj{u}{t}{\A} := \left( t, \sum_{\alpha \in \A_1} u_{(\alpha,1)} t^\alpha, \ldots, \sum_{\alpha \in \A_{n'}} u_{(\alpha,n')} t^\alpha \right).
\]
%will be used in the future.

Of course, not all coordinate hyperplanes $\A$ will give rise to restricted operators $T_\A$ which have nontrivial $L^p-L^q$ boundedness properties.  In particular, if $T_\A$ is to be bounded from any $L^p$ to some $L^q$, it must be the case that the following conditions are satisfied (the proof of necessity will be taken up in section \ref{necessitysec}):
\begin{enumerate}
\item (Dimensionality) There exists an integer $\# \A$ such that, for any $j=1,\ldots,n'$, the cardinality of $\A_j$ is $\# \A$.
\item (Scaling) There exists an integer $|\A|$ such that $ \sum_j \sum_{\alpha \in \A_j} \alpha = |\A| \mathbf 1$, where $\mathbf 1 := (1,\ldots,1)$. %(|\A|,\ldots,|\A|)$.
\item (Spanning) The multiindices $\bigcup_j \A_j$ span $\R^n$ as vectors.
\end{enumerate}
Throughout this paper, the collection $\A$ will be called {\it admissible} when it satisfies the dimensionality and scaling conditions, along with a slightly stronger form of the spanning condition (to be addressed in section \ref{vandermondesec}) and a further ``nondegeneracy'' condition:
\begin{enumerate}
\item[4] (Nondegeneracy) For each $j=1,\ldots,n'$, ${\mathbf 0} \in \A_j$, where ${\mathbf 0} = (0,\ldots, 0)$.
\end{enumerate}
With these definitions, the second main theorem of this paper is:
\begin{theorem}
Let $\A$ be admissible. Then $T_\A$ maps $L^\frac{|\A| + \# \A}{\# \A} \rightarrow L^{|\A| + \# \A}$. \label{theorema}
\end{theorem}

This theorem implies the global result of theorem \ref{maintheorem}, since the collection $\A = M_{n,d} \times \{1, \ldots, n'\}$ (which corresponds to the trivial restriction of $T_{n,n',d}$ to the whole space on which it is defined) is readily checked to be an admissible collection.

\subsection{Examples}

Theorems \ref{maintheorem} and \ref{theorema} have several corollaries which are interesting in their own right.  The first of these worth mentioning is similar to an earlier result of Ricci and Travaglini \cite{rt2001}:
\begin{corollary}
Given an $n'$-tuple of polynomials $r := (r_1,\ldots,r_{n'})$ in $n$ variables of degree at most $d$ (for $d \geq 2$), let $\mu_r$ be the measure on the graph of $r$ (as a subset of $\R^n \times \R^{n'}$) given by
\begin{equation}
\int \varphi(x,y) d \mu_r(x,y) := \int_{\R^{n}} \varphi\left(t,r_1(t),\ldots,r_{n'}(t) \right) dt.
\end{equation}
Then given any $f \in L^{p}(\R^n \times \R^{n'})$, $p:=1 + \frac{n'd}{n+1}$,
\begin{equation}
|| f \star \mu_r ||_{L^{q}(\R^n \times \R^{n'})} < \infty \label{ricci}
\end{equation}
for almost every $r$, where $q := \binom{n + d }{ d} p$.
\end{corollary}
\begin{proof}
The idea of the proof is to express integration over $y \in \R^{n'}$ and $x \in \R^n$ as integration over coefficients of polynomials and thereby reduce \eqref{ricci} to theorem \ref{maintheorem}.  Understanding $y$ is straightforward, but the $x$ variable has more subtle properties which must be understood.  The idea is that the family of polynomials $r(t+x)$ (with parameter $x$) is quite often, though not always, an $n$-dimensional hypersurface in the space of polynomials.  When this is the case, standard change-of-variables arguments may be employed.

 By definition of convolution and an elementary change of variables, one has the formula
\[ f \star \mu_r(x,y) = \int f(-t,y - \Phi_x (r)(t)) dt, \]
where $\Phi_\cdot$ is the group of transformations of polynomials given by $\Phi_x(r)(t) := r(t+x) $.  %Now for almost every $r$ (as a function of the coefficients), there is no $n$-tuple $c$ of real numbers $c_1,\ldots,c_n$ such that
%$\sum_{j=1}^n c_j \partial_j r_k(t)$ has degree less than $d-1$ for any $k=1,\ldots,n'$.  
Fix any  multiindices $\alpha_1,\ldots,\alpha_n$ of degree $d-1$ and integers $j_1, \ldots, j_n$, each in $1,\ldots,n'$, and consider the mapping from $\R^n$ to $\R^n$ given by $x \mapsto (\partial^{\alpha_1} [\Phi_x(r)]_{j_1}(0),\ldots,\partial^{\alpha_n} [\Phi_x(r)]_{j_n}(0))$ (here the derivatives are in $t$, then $t$ is set equal to zero).  This mapping is an affine linear function of $x$.  Furthermore, the entries in the Jacobian matrix each depend linearly on the coefficients of the degree $d$ terms of $r$.  %For this choice of $\alpha$'s and $j$'s, 
Standard results from algebraic geometry dictate that the corresponding mapping $x \mapsto (\partial^{\alpha_1} [\Phi_x(r)]_{j_1}(0),\ldots,\partial^{\alpha_n} [\Phi_x(r)]_{j_n}(0))$ will either be degenerate for all $r$ or invertible for almost every $r$ (comprising an open set).

Consider the choice $\alpha_1 := (d-1,0,\ldots,0), \ldots, \alpha_n := (0, 0 ,\ldots, 0, d-1)$ and $j_1 = \cdots = j_n = 1$.  For any $n'$-tuple of polynomials $r$ with $r_1(t) = \sum_{j=1}^n t_j^d$, the mapping considered in the previous paragraph is simply $x \mapsto d(x_1, \ldots, x_n)$, which is manifestly invertible.  Therefore, for this particular choice of $\alpha$'s and $j$'s, the mapping $x \mapsto (\partial^{\alpha_1} [\Phi_x(r)]_{j_1}(0),\ldots,\partial^{\alpha_n} [\Phi_x(r)]_{j_n}(0))$ must be invertible for any $r$ in some open set containing almost every $n'$-tuple of polynomials; call the set of such $r$'s the set of invertibility. The standard change-of-variables argument gives that for any sufficiently small ball $B$ centered at an $r$ in the set of invertibility, there is a constant $C$ such that, for any function $g$ depending on $n'$-tuples of polynomials,
\[ \int_B \int \int g(\Phi_x( \tilde r)) dx dy d\tilde r \leq C ||g||_1; \]
here $d\tilde{r}$ is the Lebesgue integral over the coefficients of all $n'$-tuples of polynomials in $B$. To obtain this inequality, one changes the order of integration so that the integrals in $x$ and $y$ are performed before the $\alpha_k$'th coefficient of the $j_k$-th polynomial for each $k=1,\ldots, n$ and before the constant coefficients $\alpha = \mathbf 0$ of each polynomial as well.  If the {\it remaining} coefficients (comprising $\tilde r$) are collectively referred to as $v$, the mapping $(x,y,v) \mapsto \Phi_x(r)$ has a constant, nonzero Jacobian determinant (for each value of the frozen coefficients); after the change-of-variables $(x,y,v) \mapsto \Phi_x(r)$, the integral $dx dy dv$ is simply an integral over the space of ($n'$-tuples of) polynomials.
Taking $g(r) := \left( \int f(-t,y - r(t)) dt \right)^q$ gives
\[\int_B ||f \star \mu_{\tilde r}||^q_q d \tilde r \leq C ||T_{n,n',d} f||_q^q, \]
which, by theorem \ref{maintheorem}, implies \eqref{ricci} for almost every $\tilde r$ near $r$.
\end{proof}

Another important corollary of theorem \ref{theorema} deals with restrictions $T_{\A}$ which are not overdetermined at all.  This happens, for example, when $n' = 1$ and $\# \A = n$.  In this case, theorem \ref{theorema} reduces to the following:
\begin{corollary}
Let $\alpha_1,\cdots,\alpha_n$ be multiindices on $\R^n$ which are linearly independent as vectors and sum to $(\sigma,\ldots,\sigma)$ for some integer $\sigma$.  Then the averaging operator on $\R^{n} \times \R$ given by
\[Rf(y_0,\ldots,y_{n}) := \int f(t,y_0 + y_1 t^{\alpha_1} + \cdots + y_n t^{\alpha_n}) dt \]
is bounded from $L^{\frac{n+1+\sigma}{n+1}}$ to $L^{n +1+ \sigma}$. \label{nonoverdet}
\end{corollary}

\subsection{About the proof}

As mentioned earlier, the proof is based on combinatorial tools introduced by Christ \cite{christ1998} and expanded upon by Tao and Wright \cite{tw2003} and many others.  As in these works, the idea is to consider the bilinear form induced by \eqref{theop} and iterate kernel flows of the corresponding projection operators.  Because of the concrete nature of the operators \eqref{theop}, a coordinate dependent approach will be used, and much of the general geometry found, for example, in Tao and Wright \cite{tw2003} or earlier work of the author \cite{gressman2006}, will be suppressed.

As in \cite{gressman2006}, an essential feature of the proof is that the kernel flows are ``lifted'' to a higher dimensional space in order to make the necessary change-of-variables arguments.  A new feature introduced here is that all inequalities are shown to behave well with respect to tensor products (meaning that the ``lifted'' inequalities are again lifted to product spaces of arbitrarily high dimension).  This allows one to deduce strong-type inequalities from the tensored restricted weak-type inequalities, as in earlier work of Bennett, Carbery, and Wright \cite{bcw2005} and Carbery \cite{carbery2004}.  Unlike these earlier situations, however, there is no natural ``tensor-invariance'' to exploit.  Instead, there are several new, nontrivial estimates which must be established to reproduce the earlier argument in the current, more general, case.

The rest of this paper proceeds as follows:  section \ref{measuresec} is devoted to establishing a number of inequalities concerning the Lebesgue measure (and certain $L^p$-spaces for $0 < p < 1$) which are essentially combinatorial in nature.  These inequalities will be necessary to establish the main theorem using the standard approach of geometric combinatorics.  Section \ref{tensorsec} is concerned with the introduction of tensor-product inequalities.  In particular, it is demonstrated in this section how one deduces strong-type inequalities from tensored restricted weak-type inequalities, and the tensor-product behavior of certain inequalities from section \ref{measuresec} is addressed as well.  Section \ref{sufficientsec} gives the proof of theorem \ref{theorema} (and hence theorem \ref{maintheorem} as well), and section \ref{necessitysec} establishes the necessity of various conditions of the main theorems.

\section{Measure inequalities}
\label{measuresec}
\subsection{Interpolation of monomial-weight measures}
\label{interpsec}
The first result of this section is an interpolation inequality %inequality for 
for measures on $\R^n$ which are equal to the Lebesgue measure times a monomial weight.  %{\bf MAYBE NAGEL???} To be precise, 
Let $ d \mu(x) := |x_1|^{-1} \cdots |x_n|^{-1} dx_1 \cdots dx_n$, and  for any $s \in \R^{n}_+$, let
\[ |E|_s := \int_E |x^s| d \mu(x) := \int_E \prod_{j=1}^n |x_j|^{s_j - 1} dx. \]
To prove an interpolation inequality for these measures $| \cdot |_s$, the first step is to determine the measure-theoretic properties of certain extremal sets. % which will turn out to be extremal.
%The immediate goal at this point is to establish lemma \ref{interplemma}, which gives a sharp interpolation inequality for this family of measures.  
Since it is relatively straightforward, the sharp constants are given both in proposition \ref{shapevol} and lemma \ref{interplemma}.

\begin{proposition} 
Let $v_1,\ldots,v_n \in \R^n$ and $s \in \R_+^n$ be such that $v_1,\ldots,v_n$ are linearly independent and $s$ is in the interior of the convex cone generated by $v_1,\ldots,v_n$.  Let $\sigma \in \R_+^n$ be such that $s = \sum_{j=1}^n \sigma_j v_j$.  Then for any $a \in \R_+^n$, the set $E_v^a := \set{x \in \R^n}{\sum_{i=1}^n a_i^{-1} |x|^{v_i} \leq 1}$ satisfies \label{shapevol}
\begin{equation}
|E_v^a|_s =  2^n \frac{a^\sigma}{V} \frac{\prod_{j=1}^n \Gamma(\sigma_j)}{\Gamma(1 + |\sigma|)} \label{vol1}
\end{equation}
where $V$ is the absolute value of the determinant of the matrix with columns $v_1,\ldots,v_n$.
\end{proposition}
\begin{proof}
In the integral $\int \chi_{E_v^a}(x) \prod_{j=1}^n |x_j|^{s_j-1} dx$, make the change of variables $u_i = a_i^{-1} x^{v_i}$.  The linear independence of the $v_i$'s guarantees that this map is one-to-one and onto on $\R_+^n$ (by symmetry, this is the only orthant on which the integral need be computed).  The $(i,j)$-entry of the Jacobian matrix of this change is precisely $a_i^{-1} v_{i,j} x^{v_i} x_j^{-1}$, so the absolute value of the determinant is $V \prod_{j=1}^n x_j^{-1} \prod_{i=1}^n u_i$.  It follows that $x^s d \mu(x) = a^{\sigma} V^{-1} u^\sigma d \mu(u)$ and hence
\[|E_v^a|_s = 2^n a^\sigma V^{-1} \int_T u^\sigma {d\mu(u)} \]
where $T$ is the simplex $\set{u \in \R_+^n}{\sum u_i \leq 1}$.  A straightforward induction argument (using Euler's identity for the Beta function) computes this integral and gives \eqref{vol1}.
\end{proof}

%\[\int_{E^a_v} a_j^{-1} |x|^{v_j+s} d \mu(x) = |E_v^a|_s \frac{\sigma_j}{1 + |\sigma|}\]

\begin{lemma}
Let $w_1,\ldots,w_n \in \R_+^n$ be linearly independent, and let $s \in \R_+^n$ be in the interior of the convex hull of the $w_j$'s and $0$.  Let $\theta_0,\ldots,\theta_n$ be such that $\sum_{j=0}^n \theta_j = 1$ and $s = \sum_{j=1}^n \theta_j w_j$.  Then \label{interplemma}
\begin{equation}
|E|_s   \leq   2^{n \theta_0} \prod_{j=0}^n \theta_j^{-\theta_j} \left[  \frac{\prod_{j=1}^n \Gamma(\theta_j \theta_0^{-1})}{W \Gamma(\theta_0^{-1})}\right]^{\theta_0} \prod_{j=1}^n |E|_{w_j}^{\theta_j} \label{vol2}
\end{equation}
where $W$ is the absolute value of the determinant of the matrix with columns $w_1,\ldots,w_n$.
\end{lemma}
\begin{proof}
Let $a_j := |E|_{w_j}$ and $g(x) := \sum_{j=1}^n \theta_j a_j^{-1} |x|^{w_j-s}$.  
By construction, $\int_{E} g(x) |x|^s d \mu(x) = 1 - \theta_0$.  Let $G_\lambda := \set{x \in \R^n}{g(x) < \lambda}$.  The quantity $|G_\lambda|_s$ can be evaluated using proposition \ref{shapevol}: simply take $v_j := w_j - s$ and $\sigma_j := \theta_j \theta_0^{-1}$.  Elementary computations give that $V = \theta_0 W$ (the matrix of $v$'s ``factors'' as the $w$-matrix and a matrix involving only the $\theta$'s). Combining these observations, proposition \ref{shapevol} dictates that
\begin{equation}
|G_\lambda|_s^{\theta_0} = 2^{n\theta_0} \lambda^{1- \theta_0} \prod_{j=0}^{n} \theta_j^{-\theta_j} \left[ \frac{\prod_{j=1}^n \Gamma(\theta_j \theta_0^{-1} )}{ W \Gamma(\theta_0^{-1})} \right]^{\theta_0}  \prod_{j=1}^n |E|_{w_j}^{\theta_j}. \label{interplem1}
\end{equation}  Along these same lines, $\int_{G_\lambda} g(x) |x|^s d \mu(x) = \lambda |G_\lambda|_s (1- \theta_0)$; this is because $|G|_{w_j} = \lambda |G|_s |E|_{w_j}$ by proposition \ref{shapevol} as well.
Now choose $\lambda$ so that $|G_\lambda|_s = |E|_s$.  Then $1-\theta_0 = \int_E g(x) |x|^s d \mu(x) \geq \int_{G_\lambda} g(x) |x|^s d \mu(x)$.  To see this, simply observe that the integral over $E \setminus G_\lambda$ is necessarily greater than the integral over $G_\lambda \setminus E$ (and therefore the value of the integral decreases if all the parts of $E$ outside $G_\lambda$ are moved inside $G_\lambda$).  %Invoking proposition \ref{shapevol} again, one can easily compute that 
Since $\int_{G_\lambda} g(x) |x|^s d \mu(x) = \lambda |E|_s (1 - \theta_0)$, it follows that $\lambda |E|_s \leq 1$.  Multiplying both sides of \eqref{interplem1} by $|E|_s^{1-\theta_0}$, recalling that $|G_\lambda|_s = |E|_s$ and using the inequality $\lambda |E|_s \leq 1$ gives \eqref{vol2}.
\end{proof}
In practice, the following corollary of lemma \ref{interplemma} will be more useful than lemma \ref{interplemma} itself.  At this point, accounting for constants becomes a chore and will be neglected.
\begin{corollary}
Let $w_1,\ldots,w_N \in \R_+^n$ for $N \geq n$ have a sub-$n$-tuple which is linearly independent.  Then for any positive $\theta_0,\ldots,\theta_N$ satisfying $\sum_{j=0}^N \theta_j = 1$, there exists a constant $C < \infty$ such that
\begin{equation}
|E|_s \leq C \prod_{j=1}^N |E|_{w_j}^{\theta_j} \label{interpmany}
\end{equation}
where $s = \sum_{j=1}^N \theta_j w_j$.
\end{corollary}
\begin{proof}
Suppose that $w_1,\ldots,w_n$ are linearly independent.  If $N = n$, then lemma \ref{interplemma} gives precisely the desired conclusion; assume, then, that $N > n$.  Let $\varphi := \sum_{j=n+1}^N \theta_j$.  H\"{o}lder's inequality immediately gives that
\begin{equation}
|E|_{\varphi^{-1} \sum_{j=n+1}^N \theta_j w_j} \leq \prod_{j=n+1}^N |E|_{w_j}^{\varphi^{-1} \theta_j}. \label{corr0}
\end{equation}
On the other hand, lemma \ref{interplemma} gives that
\begin{equation}
|E|_{(1 - \varphi)^{-1} \sum_{j=1}^n \theta_j w_j} \leq C \prod_{j=1}^n |E|_{w_j}^{(1-\varphi)^{-1} \theta_j} \label{corr1}
\end{equation}
because the $w$'s satisfy the independence condition, and $(1-\varphi)^{-1} \sum_{j=1}^n \theta_j = (1-\varphi)^{-1}(1 - \varphi - \theta_0) < 1$.  But H\"{o}lder's inequality also dictates that $|E|_s \leq |E|_{\varphi^{-1} \sum_{j=n+1}^N \theta_j w_j}^\varphi |E|_{(1 - \varphi)^{-1} \sum_{j=1}^n \theta_j w_j}^{1 - \varphi}$, so taking a convex combination of \eqref{corr0} and \eqref{corr1} gives \eqref{interpmany}.
\end{proof}
\subsection{Vandermonde means}
\label{vandermondesec}
Suppose $\A$ is an admissible subset of $M_{n,d} \times \{1,\ldots,n'\}$.  Let $x_j \in \R^n$ for $j=1,\ldots,\# \A$, and for $k=1,\ldots,n'$, let $V_k(x)$ be the determinant of the $\# \A \times \# \A$ matrix whose $j,m$-entry is $x_j^{\alpha_{m}}$ (where $\alpha_1,\ldots,\alpha_{\# \A}$ is an enumeration of $\A_k$).  The product of these functions $V_k$ will be called the {\it Vandermonde polynomial associated to $\A$}, and denoted $V_{\A}$, that is:
\begin{equation}
V_{\A}(x) := \prod_{k=1}^{n'} V_k(x). \label{vandermonde}
\end{equation}
If, for example, $n=n'=1$ and $\A = {0,1,\ldots,d}$, then $V_{\A}$ is exactly the $d$-th classical Vandermonde polynomial (modulo, of course, a factor of $\pm 1$).

To establish theorem \ref{theorema}, it will be necessary to have an estimate for the expectation of $|V_\A(x)|$ when the $x_j$'s ($j=1,\ldots,\# \A$) are randomly chosen points.  It will suffice for the purposes here to prove that
\begin{equation}
\int |V_\A(x)| \prod_{j=1}^{\# \A} |f_j(x_j)| dx_1 \cdots dx_{\A} \geq c \prod_{j=1}^{\# \A} ||f_j||_{L^p(\R^n)} \label{coerc}
\end{equation}
where $p := \frac{\# \A}{\# \A + |\A|} < 1$ (many more such inequalities are, in fact, true, but will not be needed here).  Along the way, the strengthened spanning condition for $\A$ will be encountered.

The proof of \eqref{coerc} begins with a definition.  Given a Lebesgue-measurable set $E \subset \R^n$, let $S^j(E)$ be the set given by %$S^j(E)$ be the Steiner symmetrization of $E$ with respect to the coordinate hyperplane $P_j := \set{(y_1,\ldots,y_d) \in \R^d}{y_j = 0}$, i.e., 
\[\set{ (y_1,\ldots,y_n) \in \R^n}{|y_j| < \frac{1}{2} \int \chi_E(y_1,\ldots,y_{j-1},s,y_{j+1},\ldots, y_n) ds}.\]
By Fubini's theorem, this set is well-defined (up to a set of measure zero); it is called the {\it Steiner symmetrization} with respect to the hyperplane $P_j := \set{(y_1,\ldots,y_n) \in \R^n}{y_j = 0}$.  Observe that the intersection of $S^j(E)$ with any line $\ell$ pointing in the $j$-th coordinate direction is simply a line segment with center in $P_j$.  Moreover, the measure of $S^j(E)$ is the same as the measure of $E$, and if $f(y)$ is any nonnegative measurable function which does not depend on $y_j$, then $\int_E f(y) dy = \int_{S^j(E)} f(y) dy$ (both of these facts follow almost directly from Fubini's theorem).  %This fact is a crucial piece 
The following propositions illustrate how Steiner symmetrization will be useful here.  Heuristically, if one wants to estimate the integral of a function $|f|$ on a set $E$, the function $f$ may be replaced by a simpler function if, in exchange, the set $E$ is replaced by a Steiner-symmetrized version of itself.
\begin{proposition}
Let $I \subset \R$ be some (possibly infinite) interval, and let $f$ be a function in $C^k(I)$ (for some fixed $k \geq 1$) which satisfies the inequality $f^{(k)}(t) \geq 1$ for all $t \in I$.  There exists a constant $c_k > 0$ such that, for any measurable set $E \subset I$, \label{sublevelprop}
\begin{equation}
 \int_E |f(t)| dt \geq c_k \int_{S(E)} \frac{|t|^k}{k!} dt. \label{steinerint}
\end{equation}
\end{proposition}
\begin{proof}
The first step is to establish the inequality \eqref{steinerint} in the case when $E$ is an interval.  Let $E := [a,b] \subset I$, $\delta := \frac{b-a}{k+1}$, and $a_j := a + j \delta$ for $j=0,\ldots,k+1$.  By the fundamental theorem of calculus (and an elementary induction argument on $k$), the following identity holds for all $f \in C^k(I)$ when $[a,b] \subset I$:
\[ \sum_{j=0}^{k} (-1)^{k-j} \binom{k}{j} \int_{a_j}^{a_{j+1}} f(t) dt = \int_0^\delta \! \! \! \cdots \! \int_0^\delta f^{(k)}(a + t_1 + \cdots + t_{k+1}) dt_1 \cdots dt_{k+1}.  
\]
The right-hand side is at least $\delta^{k+1}$ by virtue of the pointwise estimate for $f^{(k)}$; majorizing the left-hand side in the standard way gives that
\begin{equation}
  \binom{k}{j_0} \int_a^b |f(t)| dt \geq \frac{|b-a|^{k+1}}{(k+1)^{k+1}} \label{subint}
\end{equation}
where $j_0 = \frac{k}{2}$ or $\frac{k+1}{2}$ when $k$ is even or odd, respectively.  This inequality \eqref{subint} is precisely \eqref{steinerint} by virtue of the fact that $|E| = |b-a|$.

Let $F_\lambda := \set{t \in I}{|f(t)| \leq \lambda}$.  The set $F_\lambda$ is a union of no more than $k$ disjoint, closed intervals (between any two connected components of $F_\lambda$, $f'$ must vanish; but $f'$ can vanish only $k-1$ times by Rolle's theorem).  The inequality \eqref{subint} can thus be used to estimate $\int_{F_\lambda} |f(t)| dt$ by restricting to subintervals and summing.  The right-hand side of a $k$-fold sum of \eqref{subint} is minimized (given $|F_\lambda|$) when each subinterval has measure $\frac{1}{k} |F_\lambda|$; thus
\begin{equation}
 \int_{F_\lambda} |f(t)| dt \geq \frac{j_0!(k-j_0)!}{k! k^k (k+1)^{k+1}} |F_\lambda|^{k+1}. \label{subest}
\end{equation}
Now $|F_\lambda|$ is clearly a nondecreasing function of $\lambda$; as $\lambda \rightarrow \infty$, it must also be the case that $|F_\lambda| \rightarrow |I|$ (because $|f|$ is bounded on any finite subinterval of $I$).  Inner and outer regularity of the Lebesgue measure, coupled with the fact that $|\set{t \in I}{|f(t)| = \lambda}| = 0$ (solutions to $f(t) = \pm \lambda$ are isolated thanks to the derivative inequality), give that $|F_\lambda|$ is a continuous function of $\lambda$ as well.  Therefore, given an arbitrary measurable set $E \subset I$, there exists a unique minimal $\lambda > 0$ such that $|E| = |F_\lambda|$.  But
\begin{equation*} \int_E |f(t)| dt \geq \int_{E \cap F_\lambda} |f(t)| dt + \lambda |E \setminus F_\lambda|
\end{equation*}
since $|f(t)| > \lambda$ outside $F_\lambda$.  Since $|E| = |F_\lambda|$, $|E \setminus F_\lambda| = |F_\lambda \setminus E|$; then $\lambda |E \setminus F_\lambda| \geq \int_{F_\lambda \setminus E} |f(t)| dt$.  Therefore, among all measurable sets in $I$ with measure $|E|$, the integral on the left-hand side of \eqref{steinerint} is minimized for $F_\lambda$.  Furthermore, $\int_{S(E)} |t|^k dt = 2^{-k} |E|^{k+1}/(k+1)$ (by definition of the Steiner symmetrization).  Combining with \eqref{subest} gives
\[ \int_E |f(t)| dt \geq 2^k \frac{j_0! (k-j_0)!}{k^k (k+1)^k} \int_{S(E)} \frac{|t|^k}{k!} dt, \]
which is precisely the desired inequality.
\end{proof}

At this point, an ordering must be imposed on the set of multiindices.
Throughout the remainder of this section, the dictionary order on multiindices will be used; that is, 
given two multiindices $\alpha$, $\beta$, of length $n$, one says that $\alpha \leq \beta$ if and only if the smallest index $i$ for which $\alpha_i \neq \beta_i$ (if it exists) satisfies $\alpha_i < \beta_i$.  It is an elementary exercise to check that this does, in fact, define a total ordering on the set of multiindices of a given length, and consequently, any finite subset $\A_j$ of such multiindices has a maximal element.

The next proposition contains the main inequality which will be needed to establish \eqref{coerc}.  In essence, it is a generalization of proposition \ref{sublevelprop} to higher dimensions.  As before the trade-off is that integrals of (somewhat arbitrary) functions are replaced by integrals of monomials over a symmetrized set---in this case, the symmetrization is an $n$-fold symmetrization with respect to the $n$ coordinate directions (and the order of the symmetrizations is determined by the order on multiindices that was just chosen).
\begin{proposition}
For any positive integer $n'$, and for $j = 1,\ldots,n'$, let $\A_j$ be a collection of multiindices (of length $n$) such that $\# \A_j = \# \A_1 < \infty$ for all $j$. For each $j$, let $\max(\A_j)$ be the maximal element of $\A_j$, and let $\A'_j := \A_j \setminus \{ \max(\A_j) \}$.  Then there exists a constant $c$ depending only on the $\max(\A_j)$'s such that, for any measurable set $E \subset \R^n$ and any $x' \in (\R^n)^{\# \A_1 -1}$,
\[ \int_E \left| \prod_{j=1}^{n'} V_{\A_j} (x,x') \right| dx \geq c \left| \prod_{j=1}^{n'} V_{\A'_j}(x') \right| \int_{E^*} \left| x^{\sum_{j=1}^{n'} \max(\A_j)} \right| dx, \]
where $\prod_{j} V_{\A_j}$ is the Vandermonde polynomial associated to $\A$, and $E^* := S^n( \cdots (S^1(E)) \cdots)$. \label{steinerlemma}
\end{proposition}
\begin{proof}
This lemma follows from repeated application of the sublevel inequality of proposition \ref{sublevelprop}.  For any $k > [\max(\A_j)]_1$, it must be true that $\frac{\partial^k}{\partial x_1^k} V_{\A_j} (x,x') = 0$.  This is because $V_{\A_j}$ is linear in the monomials $x^\alpha$ for $\alpha \in \A_j$; if the derivative did not vanish, it would contradict the maximality of $\max(\A_j)$.  Let $\beta := \sum_{j=1}^{n'} \max(\A_j)$.  Suppose $x = (x_1,\ldots,x_n)$; by proposition \ref{sublevelprop},
\[ \int \chi_E(x) \left| \prod_{j=1}^{n'} V_{\A_j} (x,x') \right| dx_1 \geq c_{\beta_1} \int \chi_{S^1(E)}(x) \frac{|x_1^{\beta_1}|}{\beta_1!} \left| \frac{\partial^{\beta_1}}{\partial x_1^{\beta_1}} \prod_{j=1}^{n'} V_{\A_j}(x,x') \right| dx_1 \]
since the differentiated quantity on the right-hand side is independent of $x_1$.  Repeating for $x_2,\ldots,x_n$, it must be the case that
\[ \int_E  \left| \prod_{j=1}^{n'} V_{\A_j} (x,x') \right| dx \geq \prod_{i=1}^n c_{\beta_i} \int_{E^*} \frac{|x^\beta|}{\beta!} \left| \left( \frac{\partial}{\partial x} \right)^\beta  \prod_{j=1}^{n'} V_{\A_j}(x,x') \right| dx. \]
The differentiated quantity on the right-hand side is independent of $x$ entirely; moreover, by the Leibniz rule, 
\[ \left( \frac{\partial}{\partial x} \right)^\beta  \prod_{j=1}^{n'} V_{\A_j}(x,x') =  \mathop{\sum \cdots \sum}_{\gamma^1 + \cdots + \gamma^{n'} = \beta}  \left[ \prod_{j=1}^{n'} \frac{\beta!}{\gamma^j!} \left( \frac{\partial}{\partial x} \right)^{\gamma^j} V_{\A_j}(x,x') \right].\]
Since the multiindices $\gamma^j$ sum to $\beta$, and $\partial_x^{\gamma^j} V_{\A_j} \equiv 0$ if $\gamma^j > \max(\A_j)$, it follows that all terms on the right-hand side vanish except for the term $\gamma^j = \max(\A_j)$, $j=1,\ldots,n'$, and hence
\[ \int_E \left| \prod_{j=1}^{n'} V_{\A_j} (x,x') \right| dx \geq \prod_{i=1}^n c_{\beta_i} \left| \prod_{j=1}^{n'} \left(\frac{\partial}{\partial x} \right)^{\max(\A_j)} \frac{V_{\A_j}(x,x')}{(\max(\A_j))!} \right| \int_{E^*} |x^\beta| dx \]
Cramer's rule dictates that the product in absolute values on the right-hand side is precisely $|\prod_j V_{\A'_j}(x')|$, completing the proof.
\end{proof}

With the aid of proposition \ref{steinerlemma}, the desired inequality \eqref{coerc} concerning Vandermonde means is quickly established.  This is precisely the purpose of the following theorem.  In the process, the strengthened spanning condition (referred to in section \ref{mainthms}) will be employed.  Recall that the spanning condition already stated is that, for any admissible $\A$, the monomials in $\A_1,\ldots,\A_{n'}$ collectively span $\R^n$ as vectors.  The strengthened spanning condition goes as follows.  Let $\alpha_{1,k} < \alpha_{2,k} < \cdots < \alpha_{\#\A,k}$ be the ordered enumeration of $\A_k$ relative to the dictionary order, and let $\beta_j := \sum_{j=1}^{n'} \alpha_{j,k}$ for $j=1,\ldots,\# \A$.  The strengthened spanning condition holds when $\beta_1,\ldots,\beta_{\# \A}$ span $\R^n$ as vectors.  After the proof of theorem \ref{thvand}, this condition will be examined more closely, but first comes the main result of this section:
\begin{theorem}
Let $\A$ be an admissible subset of $M_{n,d} \times \{1,\ldots,n'\}$ which satisfies the strengthened spanning condition given above.  Then \eqref{coerc} holds, that is, there exists a constant $c > 0$ such that
\[ \int |V_\A(x)| \prod_{j=1}^{\# \A} |f_j(x_j)| dx_1 \cdots dx_{\# \A} \geq c \prod_{j=1}^{\# \A} ||f_j||_{L^p(\R^n)} \]
for any measurable functions $f_j$, where $p := \frac{\# \A}{\# \A + |\A|} < 1$. \label{thvand}
\end{theorem}
\begin{proof}
%Let $\beta_j := \mathbf 1 + \sum_{k=1}^{n'} \alpha_{j,k}$ for $k=1,\ldots,\# \A$.  
By repeated application of proposition \ref{steinerlemma}, it follows that there exists some $c>0$ such that, for all measurable sets $E_j \subset \R^n$, 
\[\int |V_\A(x)| \prod_{j=1}^{\# \A} \chi_{E_j}(x_j) dx \geq c \prod_{j=1}^{\# \A} |E_j^*|_{\mathbf 1 + \beta_j} \]
(using the notation of section \ref{interpsec}; here $\mathbf 1 := (1,\ldots,1)$).  Because the left-hand side is symmetric, the sets $E_2,\ldots,E_{\# \A}$ may be permuted to obtain a family of inequalities {\it while $E_1$ remains fixed}.  Taking the geometric mean of all these inequalities gives that
\[\int |V_\A(x)| \prod_{j=1}^{\# \A} \chi_{E_j}(x_j) dx \geq c |E_1^*| \left( \prod_{j'=2}^{\# \A} \prod_{j=2}^{\# \A} |E_j^*|_{\mathbf 1 + \beta_{j'}} \right)^\frac{1}{\# \A - 1} \]
(where the observation that $|E_1|_{\mathbf 1 + \beta_1} = |E_1|_{\mathbf 1} = |E_1|$ has been quietly exploited).  Since $\A$ is admissible, $ |\A| \mathbf 1 = \sum_{j=2}^{\# \A} \beta_j$, so $\sum_{j=2}^{\# \A} \mathbf 1 + \beta_j =  (\# \A - 1 + |\A|) \mathbf 1$ as well.  For any scalars $c_2,\ldots,c_{\#\A}$, \[\sum_{j=2}^{\# \A} c_j (\mathbf 1 + \beta_j) = \sum_{j=2}^{\# \A} \left(c_j + |\A|^{-1} \sum_{k=2}^{\# \A} c_k \right) \beta_j.\]
As vectors in $\R^n$, the right-hand side can assume any value for appropriate choice of the $c_j$'s (since the matrix with off-diagonal entries $|\A|^{-1}$ and diagonal entries $1 + |\A|^{-1}$ is invertible).  Therefore, the vectors $\mathbf 1 + \beta_j$ span $\R^n$ as well, and so one may employ
lemma \ref{interplemma} (or more precisely, inequality \eqref{interpmany}) to conclude that
\[\int |V_\A(x)| \prod_{j=1}^{\# \A} \chi_{E_j}(x_j) dx \geq c |E_1^*| \prod_{j=2}^{\# \A} |E_j^*|^{1 + \frac{|\A|}{\# \A -1}}
\]
(for some new constant $c$).  Recalling $|E_j^*| = |E_j|$, this inequality becomes a restricted weak-type estimate for \eqref{coerc}.  By standard machinery, this estimate may be summed to obtain Lorentz space inequalities.  In this case, the result is that
\[\int |V_\A(x)| \prod_{j=1}^{\# \A} |f_j(x_j)| dx \geq c ||f_1||_{L^1(\R^n)} \prod_{j=2}^{\# \A} ||f_j||_{L^{p_1,1}(\R^n)}
\]
where $p_1^{-1} := 1 + \frac{|\A|}{\# \A -1}$.  Now the $f_j$'s are permuted again, this time including $f_1$.  The geometric mean of these permuted inequalities is precisely
\[\int |V_\A(x)| \prod_{j=1}^{\# \A} |f_j(x_j)| dx \geq c \prod_{j=1}^{\# \A} ||f_j||_{L^1}^{\frac{1}{\# \A}}||f_j||_{L^{p_1,1}}^{1 - \frac{1}{\# \A}}.
\]
By the standard convexity inequalities for Lorentz space quasi-norms, the proof is complete.
\end{proof}
Concerning the strengthened spanning condition:  it should be noted that the strengthened spanning condition depends on the order which is imposed on multiindices.  The $\beta_j$'s themselves may change if, for example, one takes a different dictionary ordering on $M_{n,d}$ (reordering the coordinate axes, for example).  If, however, one is in the situation that $\A_j = \A_k$ for all $j,k$, the strengthened spanning condition reduces to the spanning condition mentioned in the introduction.  This is because each $\beta_j$ is necessarily equal to $n' \alpha_j$ for some $\alpha_j \in \A_j$, and changing the ordering on $M_{n,d}$ simply reorders the $\beta_j$'s.  Situations in which this occurs (and, hence, the two spanning conditions are equivalent) include that of theorem \ref{maintheorem} and corollary \ref{nonoverdet}, as well as the case of codimension $1$ averaging operators ($n'=1$).

\section{Tensor inequalities}
\label{tensorsec}
In this section, two propositions are established concerning the relationship (first observed by Carbery \cite{carbery2004}) between strong-type estimates for an operator and restricted weak-type estimates for the tensor products of that operator.  Propositions \ref{tensbound} and \ref{tenscoerc} each shed a small amount of light on this relationship from different perspectives.  The overarching idea is that weak-$L^p$ norms do not naturally behave well under tensor products, e.g., the weak-$L^p(\R^N)$ norm of $\prod_{j=1}^N f(x_j)$ is in general greater than the $N$-th power of the weak-$L^p(\R)$ norm of $f$.  If, by chance, there is some control on the growth of these norms as $N \rightarrow \infty$, then one can gain information about the regularity of $f$.  Conversely, if one has extra information about $f$ it can be possible to control the growth as $N \rightarrow \infty$.  Propositions \ref{tensbound} and \ref{tenscoerc} each establish this principle in one direction; as noted in the statements themselves, the implications go both ways.  To prove the converse of proposition \ref{tenscoerc}, one mimics the proof of proposition \ref{tensbound} and so on.  The converses have been omitted only because they are not necessary here. 
%
%The proof of proposition \ref{tensbound} was first given, in this context, by Carbery {\bf say more here!!!!!!!!!!}

To simplify notation throughout the rest of this paper, bold will be used to indicate objects in a product space.  For example, the variable $u$ will, from here on, denote a point in $\R^{\A}$.  The bold version, $\mathbf u$, is to be understood as an element of $(\R^{\A})^N$.  If an operation (which is {\it not} bold) is performed on a tensor (bold) variable, it is to be understood as a component-wise operation.  For example, 
\begin{equation*}
\proj{\mathbf u}{\mathbf t}{\A} := \left(\proj{u^1}{t^1}{\A}, \cdots, \proj{u^N}{ t^N}{\A}\right)
\end{equation*}
where, as noted, $\mathbf u := (u^1,\ldots,u^N) \in (\R^{\A})^N$ and $\mathbf t := (t^1, \ldots, t^N) \in (\R^n)^N$.  Functions which are not meant to be applied component-wise (i.e., functions which have nontrivial behavior on product spaces) and sets in product spaces will also be bold.

\begin{proposition}
The operator $T_\A$ maps $L^p(\R^{n} \times \R^{n'})$ to $L^q(\R^{\A})$ ($1 < p,q < \infty$) if (and only if) there exists a constant $C$ such that, for all positive integers $N$, \label{tensbound}
\begin{equation} 
\int_{(\R^\A)^N} {\mathbf T}_\A (\chi_{\mathbf F}) (\mathbf u) \chi_{\mathbf G}(\mathbf u) d \mathbf u \leq C^N |\mathbf F|^\frac{1}{p} |\mathbf G|^{1 - \frac{1}{q}} \label{tensboundineq}
\end{equation}
for all measurable $\mathbf F \subset (\R^n \times \R^{n'})^N$ and $\mathbf G \subset (\R^{\A})^N$, where ${\mathbf T}_\A$ is the $N$-fold tensor product of $T_\A$, i.e., 
\begin{equation*}
{\mathbf T}_\A \mathbf f(\mathbf u) = \int_{(\R^n)^N} \mathbf f \left( \proj{\mathbf u}{\mathbf t}{\A} \right) d \mathbf t.
\end{equation*}
\end{proposition}
\begin{proof}
Let $f \in L^p(\R^{n} \times \R^{n'})$ and $g \in L^{q'}(\R^{\A})$; for each integer $j$, let $F_j := \set{x \in \R^n \times \R^{n'}}{2^{j-1} \leq |f(x)| < 2^j}$ and likewise for $G_j$.  Now, for each positive integer $M$, let $f_M(x) := \sum_{|j| \leq M} 2^j \chi_{F_j}(x)$ and so on for $g_M$.  The functions $f_M$ and $g_M$ converge monotonically as $M \rightarrow \infty$ to functions majorizing $|f(x)|$ and $|g(x)|$, respectively.  Therefore (by the monotone convergence theorem)
\[\left| \int T_\A f(u) g(u) du \right| \leq \sup_{M} \int T_{\A} f_M (u) g_M(u) du. \]
For each fixed $M$ and every positive integer $N$,
\[\left(\int T_{\A} f_M (u) g_M(u) du \right)^N = \sum_{|j| \leq NM} \sum_{|j'| \leq NM} 2^{j+j'} \int {\mathbf T}_{\A} \chi_{\mathbf F_j}(\mathbf u) \chi_{\mathbf G_{j'}}(\mathbf u) d \mathbf u \]
where $\mathbf F_j := \bigcup \set{ F_{j_1} \times \cdots \times F_{j_N}}{j_1 + \cdots + j_n = j, \ |j_1|, \ldots, |j_N| \leq M}$ (and likewise for $\mathbf G_{j'}$).  To see that this is true, simply write the left-hand side as an $N$-fold product of integrals and group the terms accordingly.  By the hypothesis of this proposition, then,
\[\left(\int T_{\A} f_M (u) g_M(u) du \right)^N \leq C^N \sum_{|j| \leq NM} \sum_{|j'| \leq NM} 2^{j + j'} | \mathbf F_j|^\frac{1}{p} |\mathbf G_{j'}|^\frac{1}{q'}. \]
Applying Jensen's inequality to each sum on the right-hand side gives that the right-hand side is itself dominated by
\[ C^N (2NM+1)^{\frac{1}{p'} + \frac{1}{q}} \left( \sum_{|j| \leq NM} 2^{jp} | \mathbf F_j| \right)^\frac{1}{p} \left( \sum_{|j'| \leq NM} 2^{ j' q'} |\mathbf G_{j'}| \right)^\frac{1}{q'}. \]
Next observe that the sums over $j$ and $j'$ ({\it inside} the parentheses) are nothing other than $||f_M||_p^{Np}$ and $||g_M||_{q'}^{Nq'}$.  Taking $N$-th roots and letting $N \rightarrow \infty$, it must be the case that
\[\left| \int T_\A f(u) g(u) du \right| \leq \sup_{M} C ||f_M||_p ||g_M||_{q'}. \]
But $f_M(x) \leq 2 |f(x)|$ and $g_M(x) \leq 2 |g(x)|$, so the right-hand side is dominated by $||f||_p ||g||_{q'}$.  Therefore the operator $T_\A$ must be bounded.
\end{proof}

\begin{proposition}
Let $V$ be any $C^\infty$ function on $(\R^k)^m$, and let $0 < p_j < 1$ for $j=1,\ldots, m$.  Then there exists a constant $c > 0$ such that \label{tenscoerc}
\begin{equation}
\int | \mathbf V(\mathbf x_1,\ldots, \mathbf x_m)| \prod_{j=1}^m \chi_{\mathbf E_j}(\mathbf x_j) d \mathbf x_1 \cdots d \mathbf x_m \geq c^N \prod_{j=1}^m |\mathbf E_j|^\frac{1}{p_j} \label{tensor1}
\end{equation}
(where $| \mathbf V(\mathbf x_1,\ldots, \mathbf x_m)|$ is the component-wise product $ \prod_{j=1}^N |V(x_1^j, \ldots, x_m^j)|$) for all measurable sets $\mathbf E_j \subset (\R^k)^N$ if (and only if) there exists a constant $c > 0$ such that
\begin{equation}
\int | \mathbf V(x_1,\ldots,x_m)| \prod_{j=1}^m |f_j(x_j)| dx_1 \cdots dx_m \geq c \prod_{j=1}^m ||f_j||_{L^{p_j}(\R^k)} \label{tensor2}
\end{equation}
for all functions $f_j$ on $\R^k$.
\end{proposition}
\begin{proof}
On the left-hand side of \eqref{tensor1}, consider first the integral $dx_1^1 \cdots dx_m^1$.  By \eqref{tensor2}, there exists a constant $c$ such that
\[ \begin{split}
\int | \mathbf V(\mathbf x_1,\ldots, \mathbf x_m)| & \prod_{j=1}^m \chi_{\mathbf E_j}(\mathbf x_j) d x_1^1 \cdots d x_m^1 \\ & \geq c \prod_{i=2}^N |V(x_1^i,\ldots,x_m^i)| \prod_{j=1}^m \left( \int \left(\chi_{\mathbf E_j} (\mathbf x_j)\right)^{p_j} d x_j^1 \right)^\frac{1}{p_j}. \end{split}
\]
This inequality may be integrated with respect to $x_1^2,\ldots,x_m^2$; again \eqref{tensor2} may be applied.  Using the identity
\[ \left|\left| \left( \int \left(\chi_{\mathbf E_j} (\mathbf x_j)\right)^{p_j} d x_j^1 \right)^\frac{1}{p_j} \right|\right|_{L^{p_j}(x_j^2)} = \left( \int \left(\chi_{\mathbf E_j} (\mathbf x_j)\right)^{p_j} d x_j^1 d x_j^2 \right)^\frac{1}{p_j},\]
one proceeds by induction on $N$, arriving at \eqref{tensor1}.
\end{proof}
\section{Proof of theorem \ref{theorema}}
\label{sufficientsec}

With the machinery of the previous sections in hand, the proof of theorem \ref{theorema} may now be undertaken.
By proposition \ref{tensbound}, it suffices to show that the uniform estimate
\eqref{tensboundineq} holds.  This section is devoted to the proof of \eqref{tensboundineq}.

Throughout this section, the variable $t$ will represent a point in $\R^n$, and $u$ will represent a point in $\R^{\A}$.  For any $s \in \R^n$, let $\varphi^s_l : \R^n \times \R^{\A} \rightarrow \R^n \times \R^{\A}$ be given by
\begin{equation}
\varphi_l^{s}(t,u) := (t+s,u).
\end{equation}
Let $\A^\circ := \A \setminus \{(\mathbf 0,j) \}_{j=1,\ldots,n'}$.  For $x \in \R^{\A^\circ}$, let $\widehat{\varphi}_r^x$ be the function from $\R^n \times \R^\A$ to $\R^{\A}$ with components
\begin{equation*}
\widehat{\varphi}_r^x(t,u)_{(\alpha,j)} := \begin{cases}
u_{(\mathbf 0,j)} - \sum_{\beta \in \A_j^\circ} x_\beta t^\beta & \alpha = \mathbf 0 \\
u_{(\alpha,j)} + x_{(\alpha,j)} & \alpha \neq \mathbf 0
\end{cases}
\end{equation*}
for any $(\alpha,j) \in \A$.  Similarly, let $\varphi_r^x : \R^n \times \R^\A \rightarrow \R^n \times \R^\A$ be given by
\begin{equation}
\varphi_r^x(t,u) := (t, \widehat{\varphi}^x_r(t,u)).
\end{equation}
These maps $\varphi_l^s$ and $\varphi_r^x$ are nothing more than the kernel flow maps which appear in the work of Tao and Wright \cite{tw2003}, for example.  An important feature which distinguishes this proof from Tao and Wright's earlier work is that the flows $\varphi_l^s$ and $\varphi_l^x$ may be multidimensional flows.

Just as in the work of Christ \cite{christ1998}, one of the major components of this proof is a change-of-variables argument involving a Jacobian determinant of a repeated composition of flow maps.  In that case, the structure of the composition was fairly straightforward (arising from the repeated composition $(T^*T)^{n/2}$).  Here, in contrast, the flows are more complicated, and the change-of-variables argument takes place in a ``lifted'' space (as also occurs for X-ray like transforms \cite{gressman2006}).  
For any measurable set $F \subset \R^n \times \R^{n'}$, let
\begin{equation}
I_\A[F](t,u) := \int |V_\A(t,t+s)| \prod_{j=1}^{\# \A-1} \chi_F ( \pi_\A \varphi_l^{s_j} \varphi_r^x (t,u) ) ds dx \label{ifunc}
\end{equation} 
where $|V_\A(t,t+s)| := |V_\A(t,t+s_1,\ldots,t+s_{\#A -1})|$.  This functional is the centerpiece of the change-of-variables argument, as demonstrated by the following proposition:
\begin{proposition}
For any $(t,u) \in \R^n \times \R^{\A}$ and any measurable $F \subset \R^n \times \R^{n'} $, \label{propcov}
\[
I_\A[F](t,u) = |F|^{\# \A -1}.
\]
Furthermore, for any $(\mathbf t, \mathbf u) \in (\R^n)^N \times (\R^\A)^N$ and any measurable $\mathbf F \subset (\R^n \times \R^{n'})^N$,
\[ \mathbf I_\A[\mathbf F](\mathbf t, \mathbf u) = |\mathbf F|^{\# \A -1},
\]
where $\mathbf I_\A$ is the $N$-fold tensor product of $I_\A$.
\end{proposition}
\begin{proof}
%Fix $(t,u)$ for the moment, %and consider the map $\Phi : (\R^n)^{\# \A -1 } \times \R^{\A^\circ} \rightarrow (\R^n \times \R^{n'})^{\# \A -1}$ defined by
%\[ \Phi_j(s_1,\ldots,s_{\#A -1},x) := \pi_A \varphi_l^{s_j} \varphi_r^x(t,u) \]
%for $j=1,\ldots, \#A -1$.  Observe that the dimensionality constraint on $\A$ guarantees that the cardinality of $\A^\circ$ is just $n' (\# \A -1)$, so $\Phi$ is a map between two spaces of the same dimension.
%The first step of the proof of this proposition is to prove that the absolute value of the Jacobian determinant of $\Phi$ is precisely $|V_\A(t,t+s)|$.  
%
For each $k =1,\ldots,n'$, let $\phi^{k,s,t,u} : \R^{\A^\circ_k} \rightarrow \R^{\#A - 1}$ be the map with components 
\[\phi^{k,s,t,u}_j(x) := u_{(\mathbf 0,k)} + \sum_{\alpha \in \A^\circ_k} x_{(\alpha,k)} [(t+s_j)^\alpha - t^\alpha] + u_{(\alpha,k)} (t+s_j)^\alpha \]
for $j=1,\ldots, \# \A -1$.  The dimensionality constraint on $\A$ guarantees that $\phi^{k,s,t,u}$ is a map between spaces of the same dimension when $s,t,$ and $u$ are fixed.  The Jacobian matrix of $\phi^{k,s,t,u}$ has as its $(j,(\alpha,k))$-entry $(t+s_j)^\alpha - t^\alpha$ (that is, in the $j$-th row and the column corresponding to $(\alpha,k)$).  The absolute value of the determinant is precisely $|V_{k}(t,t+s)|$, defined at the beginning of section \ref{vandermondesec}; to see this, simply note that the former determinant can be obtained from the latter by subtracting the $t$-row from all remaining rows.  Now, in the $x$-integral appearing in \eqref{ifunc}, make the changes of variables $y_k := \phi^{k,s,t,u}(x)$.  This is permitted for almost every $s$ since $V_\A$ vanishes on a closed set of measure zero (and each $\phi^{k,s,t,u}$ depends on different $x$-variables).
Direct computation shows that
\[ \pi_A \varphi_l^{s_j} \varphi_r^x(t,u) = \left(t + s_j, \left( \phi^{k,s,t,u}_j(x)  \right)_{k = 1,\ldots, n'} \right).\]
It follows that
\[ I_\A[F](t,u) = \int \prod_{j=1}^{\# \A -1} \chi_{F}\left(t+s_j,(y_{1j},\ldots,y_{n'j})\right) ds dy. \]
After a trivial change of variables, the right-hand side is easily seen to equal $|F|^{\# \A -1}$.

As for the functional $\mathbf I_{\A}$, the exact same changes of variables can be performed on each of the factors, giving the same conclusion as in the single factor case.
\end{proof}

Before proceeding, it is useful to recall the isoperimetric formulation of restricted weak-type estimates as introduced by Tao and Wright \cite{tw2003}.
Given measurable sets $\mathbf F \subset (\R^n \times \R^{n'})^N$ and $\mathbf G \subset (\R^{\A})^N$, let $\mathbf \Omega$ be the subset of $(\R^n \times \R^{n'} \times \R^\A)^N$ given by
\begin{equation*}
\chi_{\mathbf \Omega}(\mathbf t, \mathbf u) := \chi_{\mathbf F}(\proj{\mathbf u}{\mathbf t}{\A}) \chi_{\mathbf G}(\mathbf u).
\end{equation*}
By proposition \ref{tensbound}, to prove theorem \ref{theorema}, it suffices to show that there is a constant $C < \infty$ independent of $N$, $\mathbf F$, $\mathbf G$, and $\mathbf \Omega$ such that
\begin{equation} |\mathbf \Omega| \leq C^N |\mathbf F|^\frac{\# \A}{\# \A + |\A|} |\mathbf G|^{1 - \frac{1}{\# \A + |\A|}}. \label{maindes}
\end{equation}
This inequality will be established via a careful analysis of the functional $\mathbf I_{\A}$.  By proposition \ref{propcov}, one has the identity
\[ |\mathbf F|^{\# \A -1} | \mathbf \Omega| = \int \mathbf I_{\A}[\mathbf F](\mathbf t, \mathbf u) \chi_{\mathbf \Omega}(\mathbf t, \mathbf u) d \mathbf t d \mathbf u.
\]
Using the definition of $\mathbf I_\A$ and Fubini's theorem, the order of integration can be changed so that the integration in $\mathbf t$ and $\mathbf u$ comes before the integrals in the $\mathbf s_j$'s and $\mathbf x$.
Next make the change of variables $(\mathbf t, \mathbf u) \rightarrow \varphi_r^{-\mathbf x}(\mathbf t, \mathbf u)$.  This is simply a translation of the $u$'s, so the Jacobian determinant must equal exactly $1$.% (and $\varphi^{x}_r \varphi^{-x}_r$ is the identity operator).  
Thus $|\mathbf F|^{\# \A -1} | \mathbf \Omega|$ is exactly equal to 
\begin{equation}
 \int  \left[ \int |\mathbf V_{\A}(\mathbf t, \mathbf t+ \mathbf s)| \prod_{j=1}^{\# \A -1} \chi_{\mathbf F}(\pi_\A \varphi^{\mathbf s_j}_l(\mathbf t,\mathbf u)) d \mathbf s \right] \! \! \left[ \int \chi_{\mathbf \Omega}( \varphi_r^{\mathbf x} (\mathbf t, \mathbf u)) d \mathbf x  \vphantom{\prod_{j}^{\A}} \right] d \mathbf t d \mathbf u \label{ident1}
\end{equation}
(where the integrals have again been reordered and the change $\mathbf x \rightarrow - \mathbf x$ has been made).  The following proposition will be used to estimate the second term in brackets in \eqref{ident1} so that theorem \ref{thvand} can be applied via proposition \ref{tenscoerc}:
\begin{proposition}
There exists a nonempty subset $\mathbf F' \subset \mathbf F$ such that
\begin{equation} \left[ \int \chi_{\mathbf \Omega}( \varphi_r^{\mathbf x} (\mathbf t, \mathbf u)) d \mathbf x \right] %\chi_{\mathbf G}(\mathbf u) 
\geq \frac{1}{2} \frac{|\mathbf \Omega|}{|\mathbf F|} \chi_{\mathbf F'}(\proj{\mathbf u}{\mathbf t}{\A}) \label{cut1} %\chi_{\mathbf G} (\mathbf u) 
\end{equation}
and, for any $p \geq 1$,
\begin{equation} \int \left( \int \chi_{\mathbf F'}(\proj{\mathbf u}{\mathbf t}{\A}) d \mathbf t \right)^p \chi_{\mathbf G}(\mathbf u) d \mathbf u \geq \left( \frac{1}{2} \frac{|\mathbf \Omega|}{|\mathbf G|} \right)^p | \mathbf G| . \label{cut2} 
\end{equation}
\end{proposition}
\begin{proof}
By definition of $\mathbf \Omega$ and the fact that $\pi_\A  \varphi_r^{x}(t,u) = \pi_\A(t,u)$, it must be the case that, for each pair $(\mathbf t, \mathbf u)$,
\[\int \chi_{\mathbf \Omega}( \varphi_r^{\mathbf x} (\mathbf t, \mathbf u)) d \mathbf x = \chi_{\mathbf F}(\proj{\mathbf u}{\mathbf t}{\A}) \int \chi_{\mathbf G} ( \widehat{\varphi}^{\mathbf x}_r (\mathbf t, \mathbf u)) d \mathbf x. \]
Let $\mathbf F'$ be the subset of $\mathbf F$ such that
\begin{equation} \chi_{\mathbf F'}(\proj{\mathbf u}{\mathbf t}{\A}) \int \chi_{\mathbf G} ( \widehat{\varphi}^{\mathbf x}_r (\mathbf t, \mathbf u)) d \mathbf x \geq \frac{1}{2} \frac{|\mathbf \Omega|}{|\mathbf F|} \chi_{\mathbf F'}(\proj{\mathbf u}{\mathbf t}{\A}); \label{tempeq1} \end{equation}
$\mathbf F'$ is well-defined because the integral over $\chi_{\mathbf G}$ depends only on $\pi_\A(\mathbf t, \mathbf u)$, not on $(\mathbf t, \mathbf u)$ itself.  Because $\mathbf F' \subset \mathbf F$, the left-hand side of \eqref{cut1} is greater than the left-hand side of \eqref{tempeq1}; thus \eqref{cut1} is vacuously true in this case.
On the other hand, it must also be the case that
\[ \chi_{\mathbf F \setminus \mathbf F'}(\proj{\mathbf u}{\mathbf t}{\A}) \int \chi_{\mathbf G} ( \widehat{\varphi}^{\mathbf x}_r (\mathbf t, \mathbf u)) d \mathbf x \leq \frac{1}{2} \frac{|\mathbf \Omega|}{|\mathbf F|} \chi_{\mathbf F \setminus \mathbf F'}(\proj{\mathbf u}{\mathbf t}{\A}). \]
Now integrate both sides with respect to $\mathbf t$ and $\mathbf u_{(\mathbf 0,1)}, \ldots , \mathbf u_{(\mathbf 0, n')}$ (and {\it only} these particular $\mathbf u$'s; the rest are left fixed).   A change of variables can now be made on both sides.  On the left, the change to be made is $(\mathbf t', \mathbf u') := \varphi_r^{\mathbf x}(\mathbf t, \mathbf u)$ (that is, $\mathbf t', \mathbf u'$ depend on $\mathbf x$, $\mathbf t$ and the $\mathbf u_{(0,j)}$'s); on the right, the change is $\mathbf y := \pi_\A(\mathbf t, \mathbf u)$.  Both changes are volume preserving, giving
\[ \int \chi_{\mathbf F \setminus \mathbf F'}(\proj{\mathbf u'}{\mathbf t'}{\A}) \chi_{\mathbf G} ( \mathbf u') d \mathbf t' d \mathbf u' \leq \frac{1}{2} \frac{|\mathbf \Omega|}{|\mathbf F|} |\mathbf F \setminus \mathbf F'|. \]
Subtracting both sides from $|\mathbf \Omega|$ gives \eqref{cut2} for $p = 1$.  Applying Jensen's inequality gives all remaining $p$.
\end{proof}
The proof of theorem \ref{theorema} concludes as follows:  combining the previous proposition with the identity \eqref{ident1}, it follows that $|\mathbf F|^{\# \A -1}  |\mathbf \Omega|$ is greater than or equal to 
\[\frac{1}{2} \frac{|\mathbf \Omega|}{|\mathbf F|} \int \left[ \int |\mathbf V_{\A}(\mathbf t,\mathbf t + \mathbf s)| \chi_{\mathbf F'}(\proj{\mathbf u}{\mathbf t}{\A}) \prod_{j=1}^{\# \A -1} \chi_{\mathbf F}(\pi_\A \varphi_l^{\mathbf s_j}(\mathbf t, \mathbf u)) d \mathbf s d \mathbf t \right] \chi_{\mathbf G}(\mathbf u) d \mathbf u. \]
Observe that $\pi_\A \varphi_l^{\mathbf s_j}(\mathbf t, \mathbf u) = \pi_\A(\mathbf t + \mathbf s_j, \mathbf u)$;  thus by theorem \ref{thvand} and proposition \ref{tenscoerc}, there is a constant $c > 0$ which is independent of $\mathbf F$, $\mathbf G$, $\mathbf \Omega$, and $N$ such that the quantity in brackets is at least
\[ c^N \left( \int \chi_{\mathbf F'}(\pi_\A(\mathbf t, \mathbf u) d \mathbf t \right)^{\# \A + |\A|} \]
for each value of $\mathbf u$ (since $\mathbf F$ may be everywhere replaced by $\mathbf F'$ with impunity).  Now by \eqref{cut2}, it follows that
\[ |\mathbf F|^{\# \A -1}  |\mathbf \Omega| \geq \frac{c^N}{2} \frac{|\mathbf \Omega|}{|\mathbf F|} \left( \frac{1}{2} \frac{|\mathbf \Omega|}{|\mathbf G|} \right)^{\# \A + |\A|} |\mathbf G|, \]
which is precisely of the desired form \eqref{maindes} after elementary manipulations.

\section{Necessity}
\label{necessitysec}
%{\bf NECESSITY in the main theorem}
%{\bf MORE about spanning condition}

\subsection{Admissibility criteria}
\label{admissibility}
An important condition of theorem \ref{theorema} is that the set $\A$ be admissible.  As mentioned in the introduction, many of the admissibility criteria are, in fact, necessary for $L^p$ boundedness to hold in any form at all. By scaling, it is fairly straightforward to see that the dimensionality and scaling assumptions are necessary, and that only one global $L^p \rightarrow L^q$ estimate can hold.  Let $\delta \in \R^n_+$ and $\delta' \in \R^{n'}_+$.
For any function $f$ on $\R^n \times \R^{n'}$, the $(\delta,\delta')$ dilation of $f$ is defined to be
\[f_{\delta,\delta'}(s,w) := f(t_1 \delta_1,\ldots,t_n \delta_n, w_1 \delta'_1,\ldots, w_{n'} \delta'_{n'}).\]
Likewise, for any function $g$ on the parameter space, let
\[g^{\delta,\delta'}(u) := \delta^{- \mathbf 1} g \left( \left\{ \delta'_j \delta^{-\alpha} u_{(\alpha,j)} \right\}_{(\alpha,j) \in \A} \right),\]
where $\mathbf 1$ is the multiindex $(1,\ldots,1)$.
The standard change-of-variables argument shows that $T_A(f_{\delta,\delta'}) = (T_\A f)^{\delta,\delta'}$.  Furthermore, $||f_{\delta,\delta'}||_p = \delta^{-\frac{1}{p} \mathbf 1} {\delta'}^{-\frac{1}{p} \mathbf 1} ||f||_p$ and $||g^{\delta,\delta'}||_q = \delta^{- \mathbf 1 + \frac{1}{q} v} {\delta'}^{- \frac{1}{q} v'} ||g||_q$, where $v := \sum_{(\alpha,j) \in \A} \alpha$ and $v' := (v'_1,\ldots,v'_{n'})$ with $v_j$ equal to the cardinality of $\A_j$.  By the usual arguments, for any $L^p \rightarrow L^q$ estimate to hold, it must be the case that
\[ \frac{v}{q} = \left(1 - \frac{1}{p} \right) \mathbf 1 \mbox{ and } \frac{v'}{q} = \frac{\mathbf 1}{p}. \]
Thus the scaling and dimensionality conditions on $\A$ are necessary for $T_\A$ to be bounded at all, and when satisfied, $T_\A$ can map $L^p \rightarrow L^q$ only when $p = \frac{|\A| + \# \A}{\# \A}$ and $q = |\A| + \# \A$.

Next, consider what happens when $\A$ fails to satisfy the (weak) spanning condition; that is, suppose that the monomials in $A := \bigcup_j \A_j$ span only some subspace of $\R^{n}$ (when interpreted as vectors).  %Let $E \subset \R^n$ be the set on which $1 \leq \sum_{\alpha \in A} |x^\alpha| \leq 2$, and let $E_R$ be that subset of $E$ which is contained in the ball of radius $R$ centered at the origin.
Let $\beta_1,\ldots,\beta_m$ be linearly independent monomials which span the same subspace as $\bigcup_j \A_j$, and let $\beta_{m+1},\ldots,\beta_{n}$ be linearly independent vectors such that $\beta_1,\ldots,\beta_n$ span $\R^n$.  Now let $E_R \subset \R^n$ be the set on which $1 \leq |x^{\beta_j}| \leq 2$ for $j=1,\ldots, m$ and $1 \leq |x^{\beta_j}| \leq R$ for $j=m+1,\ldots,n$.  To compute the measure of this set use the change-of-variables $y_j := x^{\beta_j}$,
as in proposition \ref{shapevol}. 
One obtains $|E_R| = C |\ln R|^{n-m}$.  Now let $f_R$ be the characteristic function on $\R^n \times \R^{n'}$ of the set $E_R$ times the ball of radius $1$ centered at the origin in $\R^{n'}$.  For all $u$ sufficiently near zero, $T_\A f_R(u)$ also grows like $|\ln R|^{n-m}$ (since $\sum u_{(\alpha,j)} t^{\alpha}$ will be bounded for all $t \in E_R$).  If it is to be the case that $||T_\A f_R||_q \leq C ||f_R||_p$, taking $R \rightarrow \infty$ shows that $p$ must be less than $q$.  This, however, cannot happen because of the dimensionality and scaling conditions.

\subsection{Local $L^p$ estimates of theorem \ref{maintheorem}}

To conclude, consider the necessity claims of theorem \ref{maintheorem}.  The necessary constraints on global boundedness of $T_{n,n',d}$ are easily established by the same scaling argument used to demonstrate the necessity of the dimensionality and scaling conditions.  The only new feature of theorem \ref{maintheorem} not present in theorem \ref{theorema} is the claim concerning local $L^p$ estimates.

For each $\delta > 0$, let $F_\delta \subset \R^n \times \R^{n'}$ and $G_\delta \subset (\R^{n'})^{M_{n,d}}$ be given by
$F_\delta := \set{(t,s)}{|t_j| \leq \delta, |s_k| \leq C \delta^l}$ and
$G_\delta := \set{u}{|u_{(\alpha,j)}| \leq \delta^{\max \{l - |\alpha|,0 \}}}$ (for some fixed constant $C$).  Elementary counting shows that $|F_\delta| = C^{n'} \delta^{n + n'l}$ and $|G_\delta| = \delta^{K}$, where $K := n' \binom{n+l}{n+1}$.  If $C$ is fixed suitably large, it will be the case that
$\int_{G_\delta} T_{n,n',d} \chi_{F_\delta} = \delta^n |G_\delta|$ since $|\sum_{|\alpha| \leq d} u_{(\alpha,j)} t^\alpha| \leq C \delta^l$ for all $u \in G_\delta$.  If $T$ is to be bounded from $L^p$ to $L^q$, it must therefore be the case that $\delta^n |G_\delta| \leq C' |F_\delta|^{1/p} |G_\delta|^{1/q'}$.  Letting $\delta \rightarrow 0$, the inequality can hold for some $C'$ only when the exponent of $\delta$ on the left-hand side is greater than the exponent on the right.  This gives the necessary inequalities
\begin{equation}
 n + \frac{ n'}{q} \binom{n+l}{n+1} \geq \frac{n + l n'}{p} \label{constr1}
\end{equation}
for $l=1,\ldots,d$.

Finally, let $F'_\delta := \set{(t,s)}{|t_j| \leq 1, |s_k| \leq C \delta}$
and $G'_\delta := \set{u}{|u_{(\alpha,j)}| \leq \delta}$.  In this case, $|F'_\delta| = C^{n'} \delta^{n'}$ and $|G'_\delta| = \delta^{K'}$ where $K' := n' \binom{n+d}{d}$.  Proceeding as before, it follows that $\int_{G'_\delta} T_{n,n',d} \chi_{F'_\delta} = |G'_\delta|$, and after computing exponents, that
\begin{equation}
\frac{n'}{q} \binom{n+d}{d} \geq \frac{n'}{p} \label{constr2}
\end{equation}
The constraints \eqref{constr1} combined with \eqref{constr2} give precisely the necessity conditions of theorem \ref{maintheorem}, illustrated in figure \ref{thefig}.

\section{Acknowledgments}
The author would like to thank E. M. Stein for his helpful comments and suggestions while preparing this paper.

%
%\[\chi_{\mathbf \Omega}(\mathbf s, \mathbf w) =  \chi_{\mathbf F} (\varphi_1^{\mathbf s}(\mathbf w)) \chi_{\mathbf G} (\mathbf w)\] % d \mathbf s d \mathbf w \]
%\[\left( \int \chi_{\mathbf \Omega}(\mathbf s, \varphi_2^{\mathbf x} (\mathbf w)) d \mathbf x \right) \chi_{\mathbf G}(\mathbf w) \geq \frac{1}{2} \frac{|\mathbf \Omega|}{|\mathbf F|} \chi_{\mathbf F^*} (\varphi_1^{\mathbf s}(\mathbf w)) \chi_{\mathbf G} (\mathbf w) \]
%
%\[
%2 |\mathbf F|^{\# \A} \geq  \int  \left( \int \! \! \! \int |V_{\A}(\mathbf s, \mathbf t)| \chi_{\mathbf F^*}(\varphi_1^{\mathbf s}(\mathbf w)) \prod_{j=1}^{\# \A -1} \chi_{\mathbf F}(\varphi_1^{\mathbf t_j}(\mathbf w)) d \mathbf s d \mathbf t \right) \chi_{\mathbf G}(\mathbf w) d \mathbf w 
%\]
%
%\[
%\begin{split}
%\int \! \! \! \int & |V_{\A}(\mathbf s, \mathbf t)| \chi_{\mathbf F^*}(\varphi_1^{\mathbf s}(\mathbf w)) \prod_{j=1}^{\# \A -1} \chi_{\mathbf F}(\varphi_1^{\mathbf t_j}(\mathbf w)) d \mathbf s d \mathbf t \\ & \geq c^{N} \left( \int \chi_{\mathbf F^*}(\varphi^{\mathbf s}(\mathbf w)) d \mathbf s \right )^{1 + \frac{|\A|}{\# \A}} \left( \int \chi_{\mathbf F}(\varphi^{\mathbf s}(\mathbf w)) d \mathbf s \right )^{\left(1 + \frac{|\A|}{\# \A}\right)(\# \A -1)} 
%\end{split}
%\]
%Remarks:  In proposition \ref{steinerlemma}, could have put the $V_{\A_j}$'s to arbitrary positive powers (the argument needs to be souped up a bit).  Also can make sense of $\# \A_j \neq \# \A_k$.
\bibliography{mybib}
\end{document}

%% file: header.tex
\usepackage{amsmath}
\usepackage{amsthm}
\usepackage{amsfonts}
\usepackage{amsbsy}
\usepackage{amssymb}
\usepackage{pstricks}
\newtheorem{theorem}{Theorem}
\newtheorem{proposition}{Proposition}
\newtheorem{corollary}{Corollary}

\newtheorem{lemma}{Lemma}

\newcommand{\R}{{\mathbb R}}

\newcommand{\set}[2]{ \left\{ #1 \ \left| \ #2 \right. \right\} }

%% file: hyprsfc.bbl
\providecommand{\bysame}{\leavevmode\hbox to3em{\hrulefill}\thinspace}
\providecommand{\MR}{\relax\ifhmode\unskip\space\fi MR }
\providecommand{\MRhref}[2]{%
  \href{http://www.ams.org/mathscinet-getitem?mr=#1}{#2}
}
\providecommand{\href}[2]{#2}
\begin{thebibliography}{10}

\bibitem{bcw2005}
Jonathan Bennett, Anthony Carbery, and James Wright, \emph{A non-linear
  generalization of the {L}oomis-{W}hitney inequality and applications}, Math.
  Rest. Lett. \textbf{12} (2005), 443--457.

\bibitem{bgt2005}
Luca Brandolini, Allan Greenleaf, and Giancarlo Travaglini,
  \emph{${L}^p-{L}^{p'}$ estimates for overdetermined {R}adon transforms},
  Transactions of the A.M.S., to appear.

\bibitem{carbery2004}
Anthony Carbery, \emph{A multilinear generalisation of the {C}auchy-{S}chwarz
  inequality}, Proc. Amer. Math. Soc. \textbf{132} (2004), no.~11, 3141--3152
  (electronic).

\bibitem{christ1984}
Michael Christ, \emph{Estimates for the {$k$}-plane transform}, Indiana Univ.
  Math. J. \textbf{33} (1984), no.~6, 891--910.

\bibitem{christ1998}
\bysame, \emph{Convolution, curvature, and combinatorics: a case study},
  Internat. Math. Res. Notices (1998), no.~19, 1033--1048.

\bibitem{drury1983}
S.~W. Drury, \emph{{$L\sp{p}$} estimates for the {X}-ray transform}, Illinois
  J. Math. \textbf{27} (1983), no.~1, 125--129.

\bibitem{drury1984}
\bysame, \emph{Generalizations of {R}iesz potentials and {$L\sp{p}$} estimates
  for certain {$k$}-plane transforms}, Illinois J. Math. \textbf{28} (1984),
  no.~3, 495--512.

\bibitem{drury1986}
\bysame, \emph{An endpoint estimate for certain {$k$}-plane transforms}, Canad.
  Math. Bull. \textbf{29} (1986), no.~1, 96--101.

\bibitem{gressman2006}
Philip Gressman, \emph{${L}^p$-improving properties of {X}-ray like
  transforms}, Math. Res. Lett. \textbf{13} (2006), no.~5, 787--803.

\bibitem{john1938}
Fritz John, \emph{The ultrahyperbolic differential equation with four
  independent variables}, Duke Math. J. \textbf{4} (1938), no.~2, 300--322.

\bibitem{os1982}
D.~M. Oberlin and E.~M. Stein, \emph{Mapping properties of the {R}adon
  transform}, Indiana Univ. Math. J. \textbf{31} (1982), no.~5, 641--650.

\bibitem{rt2001}
Fulvio Ricci and Giancarlo Travaglini, \emph{Convex curves, {R}adon transforms
  and convolution operators defined by singular measures}, Proc. Amer. Math.
  Soc. \textbf{129} (2001), no.~6, 1739--1744.

\bibitem{tw2003}
Terence Tao and James Wright, \emph{{${L}\sp p$} improving bounds for averages
  along curves}, J. Amer. Math. Soc. \textbf{16} (2003), no.~3, 605--638.

\end{thebibliography}
